\documentclass[12pt]{amsart}
\usepackage{fullpage, amsfonts, amssymb, amsmath}
\usepackage{graphicx}
\usepackage{comment}
\usepackage{blkarray}
\usepackage{appendix}
\usepackage{tikz,tkz-graph}
\usepackage{tkz-berge}

\newtheorem{thm}{Theorem}

\newtheorem{prp}[thm]{Proposition}

\newcommand{\del}{\backslash}
\newcommand{\con}{/}

\numberwithin{thm}{section}
\numberwithin{equation}{section}

\begin{document}
\title{Intersection of Longest Cycle and Largest Bond in 3-Connected Graphs}

\author{Emily Ren}
\address{Diamond Bar High School \\
Diamond Bar, California, USA}
\email{emilyrenla@gmail.com}

\begin{abstract}
A bond in a graph is a minimal nonempty edge-cut. A connected graph $G$ is dual Hamiltonian if the vertex set can be partitioned into two subsets $X$ and $Y$ such that the subgraphs induced by $X$ and $Y$ are both trees. Research in cycles and related topics is a fundamental area in graph theory, and there is much interest in studying the longest cycles and largest bonds in graphs. H. Wu conjectured that any longest cycle must meet any largest bond in a simple 3-connected graph. In this paper, the author proves that the above conjecture is true for certain classes of 3-connected graphs: Let $G$ be a simple 3-connected graph with $n$ vertices and $m$ edges. Suppose $c(G)$ is the size of a longest cycle, and $c^*(G)$ is the size of a largest bond. Then each longest cycle meets each largest bond if either $c(G)\ge n - 3$ or $c^*(G)\ge m - n - 1$. Sanford \cite{sanford2006} determined in her Ph.D. thesis the cycle spectrum of the well-known generalized Petersen graph $P(n,2)$($n$ is odd) and $P(n,3)$($n$ is even). Flynn  \cite{flynn} proved in her honors thesis that any generalized Petersen graph $P(n,k)$ is dual Hamiltonian. The author studies the bond spectrum (called the co-spectrum) of the generalized Petersen graphs and extends Flynn’s result by proving that in any generalized Petersen graph $P(n, k)$, $1\le k\ < \frac{n}{2}$, the co-spectrum of $P(n,k)$ is \{3, 4, 5, …, $n$ + 2\}. \vspace{0.2in}


\end{abstract}

\maketitle

\section{Introduction}

All graphs in this paper are simple. We follow the notation from Douglas West's textbook \cite{west2001}. A graph is {\it connected} if there is a path between every given pair of vertices. A graph is disconnected if it is not connected. An edge $e$ of graph $G$ is called {\it cut-edge} if $G - e$ has more components than $G$. Let $A \subseteq V(G)$, define $G[A]$ as the subgraph induced by $A$. A {\it bond} in a graph is a minimal nonempty edge-cut. For a connected graph $G$, the following is another way to view a bond. Given a connected graph $G$, suppose that the vertex set $V(G)$ is partitioned into two sets $X$ and $Y$. If the subgraphs induced by $X$ and $Y$, denoted as $G[X]$ and $G[Y]$, are both connected, then the set of edges between $X$ and $Y$ forms a bond, denoted as $[X, Y]$. A graph is called {\it k-connected} if it has at least $k$ + 1 vertices and remains connected whenever fewer than $k$ vertices are removed. In a graph $G$, we use $c(G)$ to denote the size of a longest cycle (also commonly referred to as the {\it circumference} of $G$) and $c^*(G)$ to denote the size of a largest bond.  The cycle spectrum of a graph is the set of all cycle lengths, and the co-spectrum of the graph is the set of all bond sizes. A graph is {\it Hamiltonian} if there is a cycle, which passes every vertex of the graph exactly once.  A graph $G$ is called a {\it dual Hamiltonian} if it contains a bond of size $|E(G)| - |V(G)|$ + 2. P. A cycle $C$ {\it meets} a bond $B$ if they have at least one common edge (indeed, they will have at least two common edges if they do meet).

 Longest cycle and largest bond problems are important in graph theory and there is extensive literature in studying these problems (see, for example, \cite{bondy}).   Indeed, research in cycles in graphs is a fundamental research area in graph theory and has many applications.  For example, the famous traveling salesman problem has many real world applications (see, for example, \cite{west2001}).  P. Wu \cite{p-wu} proved the following interesting result relating bonds and longest cycles.

\begin{thm}(P. Wu, 1997)\cite{p-wu}
\label{p-wu} 
 Let $G$ be a $2$-connected graph with circumference $c(G)$, then there is a bond meeting every cycle of length at least c(G)$-1$.
\end{thm}

Extending the above result, McGuinness \cite{McGuinness} proved the following result for $k$-connected graphs. 

\begin{thm}(S. McGuinness, 2005)\cite{McGuinness} 
\label{McGuinness} 
Let $G$ be a k-connected graph with longest cycle of length $c\ge 2k$ and $k\ge 2$. Then there is a bond meeting every cycle of length at least c $-$ k $-$ $2$.
\end{thm}

In a very recent result,  Zhao, Wei, and Wu \cite{Zhao} proved a path version of the above results for $k$-connected graphs. We state their theorem for 3-connected graphs only.

\begin{thm}(Q. Zhao, B. Wei, and H. Wu)\cite{Zhao} 
\label{WWZ} 
Let G be a $3$-connected graph of order $n\ge 6$ with longest path of length p. Then there is a bond meeting all paths of length at least p $-$ $2$.
\end{thm}

The above results show that there are some interesting links between bonds and  longest cycles in  $k$-connected graphs for $k\ge 2$. 
 Inspired by Theorem \ref{p-wu}, Haidong Wu made the following very interesting conjecture \cite{wu} several years ago. \\

\noindent \textbf{Conjecture}:  {\it In a simple $3$-connected graph, every longest cycle must meet any largest bond. }

Any bond in a graph induces a partition of the vertex set.  If true, this conjecture provides a very interesting link between longest cycles and largest bonds in 3-connected graphs. It implies that  for any largest bond, every longest cycle must sit in exactly one side of the partitioned 
sets associated with the largest bond. In this paper, we study this conjecture and prove the following result which partially proves the conjecture. 
\begin{thm}
\label{main} 
Let $G$ be a $3$-connected graph with $n$ vertices and $m$ edges. Suppose $c(G)$ is the size of a longest cycle, and $c^*(G)$ is the size of a largest bond. Then each longest cycle meets each largest bond if one of the following is true: 
\\i) $c(G)\ge n - 3$, or
\\ii) $c^*(G)\ge m - n - 1.$
\end{thm}

There is much interest in studying the Hamiltonian and dual Hamiltonian graphs in the literature. Moreover, some researchers have studied the longest cycles in 3-connected graphs \cite{bondy}.  In the Ph.D. dissertation \cite{sanford2006} by Alice Jewel Sanford in June 2006, the author determined the cycle spectrum of $P(m, 2)$ for every odd $m$. She also found the cycle spectrum of $P(m, 3)$ for $m$ relatively prime to 3, among other results. In the honors thesis “The Largest Bond in 3-Connected Graphs” \cite{flynn} by Melissa Flynn in 2017, the author proved the next theorem, which has been generalized by Costalonga \cite{costalonga}.

\begin{thm}(Flynn, 2017) \cite{flynn} 
\label{flynn} 
 The generalized Petersen graph $P(n, k)(1\le k\ <\frac{n}{2})$ is dual Hamiltonian. That is, it has a bond of the largest possible size.
 \end{thm}

In this research, we extend the above result and determine the co-spectrum of all generalized Petersen graphs. 

\begin{thm}
\label{main2} 
In any generalized Petersen graph $P(n, k)$, where $1\le k\ <\frac{n}{2}$, the co-spectrum of P(n,k) is $\{3, 4, 5……, n + 2\}$.

\end{thm}

In Section 2, we prove Theorem \ref{main2}. In Section 3, we will prove our main result Theorem \ref{main}.

\section{Co-spectrum of the Generalized Peterson Graphs}

Recall that a graph with $n$ vertices and $m$ edges is dual Hamiltonian if it has a bond of size $m-n+2$.

\begin{center} 
\includegraphics[width=110mm]{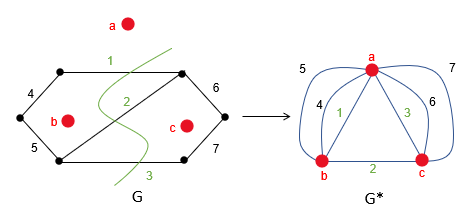}\\
Figure 1: A dual Hamiltonian graph
\end{center}

In Figure 1, a Hamiltonian graph $G$ has a bond size of $|E(G)| - |V(G)|$ + 2 = 7 $-$ 6 + 2 = 3. \{1, 2, 3\} is exactly a bond size of $|E(G)| - |V(G)|$ + 2. $G^*$ is the dual graph of $G$. In $G^*$, \{1, 2, 3\} forms a Hamiltonian cycle. Thus $G$ is a dual Hamiltonian graph. 

 The next results are well-known \cite{west2001}. 

\begin{prp}
\label{upperbound}
 Let $G$ be a connected graph with $n$ vertices and $m$ edges, then $c^*(G)\le m - n + 2$.
\end{prp}

\begin{prp}
\label{dual-H} 
Let $G$ be a connected graph with $n$ vertices and $m$ edges, then $G$ is dual-Hamiltomian (that is, $c^*(G)$ = $|E(G)| - |V(G)|$ + $2$) if and only if $V(G)$ can be partitioned as $A\cup B$, where $G[A], G[B]$ are both trees.
\end{prp}

The graph below is the famous Petersen graph, a very important graph in graph theory.

\begin{center}
\begin{tikzpicture}[rotate=90]
  \tikzset{VertexStyle/.style = {shape=circle,fill=white,
                                 minimum size = 12pt,draw}}
 \SetVertexNoLabel
 \grGeneralizedPetersen[prefix=u,prefixx=v,RA=2.6,RB=1.8]{5}{2}
 \AssignVertexLabel{u}{$x_1$,$x_2$,$x_3$,$x_4$,$x_5$}
\AssignVertexLabel{v}{$y_{1}$,$y_2$,$y_3$,$y_4$,$y_5$}
\end{tikzpicture}\\
Figure 2: A Petersen graph $P(5,2)$
\end{center}

The generalized Petersen graphs, introduced in 1950 by H.S.M. Coxeter and named in 1969 by Mark Watkins [8], are a family of cubic graphs formed by connecting the vertices of a regular polygon to the corresponding vertices of a star polygon.

The generalized Petersen graphs are a class of graphs which follow a similar design. The generalized Petersen graph, denoted $P(n, k)$ for $n\ge 3$ and $1\le k\ < \frac{n}{2}$, is defined as follows:

\begin{itemize}
\item The vertex set is $\{x_1, x_2, ..., x_n, y_1, y_2, ..., y_n\}$.
\item The edge set is composed of $\{x_ix_{i+1}, x_iy_i, y_iy_{i+k}: i = 1, ..., n\}$ where subscripts are read modulo $n$.
\end{itemize}

\begin{center} 
\includegraphics[width=80mm]{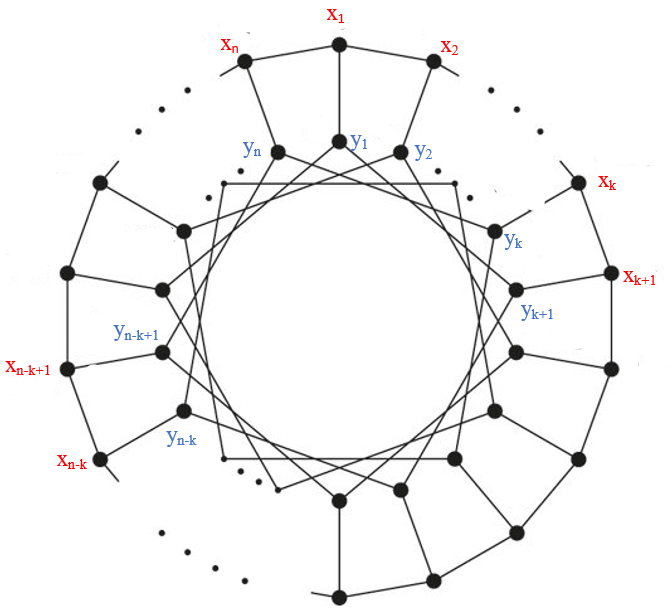}\\
Figure 3: A generalized Petersen graph $P(n, k)$
\end{center}

The generalized Petersen graph satisfies the following conditions:
\begin{itemize}
\item The graph has $2n$ vertices.
\item The graph has $3n$ edges.
\item The graph is 3-regular.
\item $P(n, k)$ is dual Hamiltonian with a largest bond $(n + 2)$, which was proved by Flynn \cite{flynn}
\end{itemize}

For example, in $P(8, 3)$ as shown in Figure 4, there are 16 vertices and 24 edges. The largest bond is 10. 

\begin{center}
\begin{tikzpicture}[rotate=90]
  \tikzset{VertexStyle/.style = {shape=circle,fill=white,
                                 minimum size = 12pt,draw}}
 \SetVertexNoLabel
 \grGeneralizedPetersen[prefix=u,prefixx=v,RA=2.8,RB=2.1]{8}{3}
 \AssignVertexLabel{u}{$x_1$,$x_2$,$x_3$,$x_4$,$x_5$, $x_6$,$x_7$,$x_8$}
\AssignVertexLabel{v}{$y_{1}$,$y_2$,$y_3$,$y_4$,$y_5$, $y_7$,$y_7$,$y_8$}
\end{tikzpicture}\\
Figure 4: A Petersen graph $P(8,3)$
\end{center}

{\bf Proof of Theorem \ref{main2}:}

\textbf{Step 1}: In Figure 5, we partition $V(P(n, k))$ into $A\cup B$. $A = \{x_1, x_2, ..., x_t\}$, $B = V(P(n, k)) - A$, where $1\le t\le n - k$. Clearly, $G[A]$ is a path (where $t$ = 1, $G[A]$ is an isolated vertex, but it is still connected), thus connected.

\begin{center} 
\includegraphics[width=80mm]{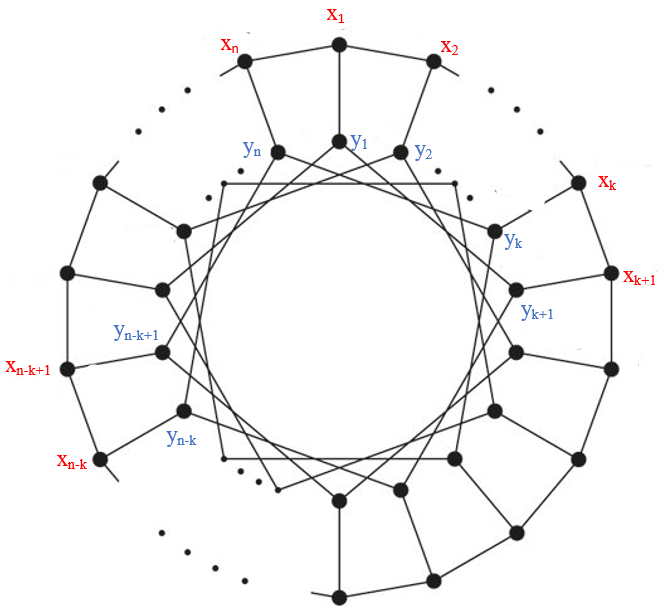}\\
Figure 5: A generalized Petersen graph $P(n, k)$
\end{center}

Let $S = \{x_{t+1}, x_{t+2}, ..., x_n\}$. Then clearly $G[S]$ is also a path. Note that each vertex $y_i (1\le i\le n)$ is connected to a vertex of $S$ by some paths in $G[B]$. Thus $G[B]$ is also connected.  Therefore, $[A,B]$ is a bond in $P(n, k)$, and this bond has size $|[A,B]| = t + 2$. Note that $t$ can be any integer from 1 to $n - k$, thus $P(n, k)$ has bonds of sizes 3, 4, 5, …, $n - k$ + 2, respectively.

\textbf{Step 2}: We show $P(n, k)$ has bonds of sizes $n - k$ + 2 + $i$ for any $i$, $1\le i\le k$.

Let $A = \{x_1, x_2, ..., x_{n-k}, y_1,  y_2, ..., y_i\}$, where $1\le i\le  k$, and $B = V(P(n, k)) - A$. We show that $[A, B]$ is a bond. Clearly $G[A]$ is connected. We show that $G[B]$ is also connected. Let $S = \{x_{n-k+1}, x_{n-k+2}, …, x_n\}$, then $G[S]$ is a path. 

For each vertex $y_j (n - k + 1\le j\le  n)$, $x_jy_j\in E(G[B])$. For each $y_j$, where $i + 1\le j\le n - k$, there is a path from $y_j$ to a vertex in $S$ in $G[B]$ by the definition of $P(n, k)(1\le k\ < \frac{n}{2})$. Hence, again $G[B]$ is connected, thus $[A, B]$ is a bond. This bond has size $|[A,B]|$ = $2i + (n -k - i) + 2 = n - k + i + 2$. Therefore, $P(n, k)$ has bonds of sizes $n - k + 3, n - k + 4, ..., n - k + k + 2 = n + 2$, respectively. 

By Step 1 and Step 2, we conclude that the co-spectrum of $P(n, k)$ is $\{3, 4, 5, …, n + 2\}$. This completes the proof of Theorem \ref{main2}.

Let's take $P(20, 4)$ as an example in Figure 6. The size of the largest bond is 22, thus it has a bond of all sizes from 3 to 22 by our theorem.

\begin{center}
\begin{tikzpicture}[rotate=90]
  \tikzset{VertexStyle/.style = {shape=circle,fill=white, minimum size = 14pt,draw}}
 \SetVertexNoLabel
 \grGeneralizedPetersen[prefix=u,prefixx=v,RA=3,RB=2.2]{20}{4}
\AssignVertexLabel{u}{$x_1$,$x_2$,$x_3$,$x_4$,$x_5$, $x_6$,$x_7$,$x_8$, $x_9$,$x_{10}$,$x_{11}$,$x_{12}$, $x_{13}$,$x_{14}$,$x_{15}$,$x_{16}$,$x_{17}$,$x_{18}$,$x_{19}$, $x_{20}$}
\AssignVertexLabel{v}{$y_{1}$,$y_2$,$y_3$,$y_4$,$y_5$, $y_6$,$y_7$,$y_8$, $y_9$,$y_{10}$,$y_{11}$,$y_{12}$, $y_{13}$,$y_{14}$,$y_{15}$,$y_{16}$,$y_{17}$,$y_{18}$,$y_{19}$, $y_{20}$}
\end{tikzpicture}\\
Figure 6: A Petersen graph $P(20,4)$
\end{center}

\section{Proof of the Main Theorem}

In this section, we prove our main result Theorem \ref{main}. We will use the following result of Dirac \cite{west2001}. 

\begin{thm}
\label{dirc} 
Let $G$ be a $2$-connected graph with minimum degree of $\delta$, then $c(G)\ge min\{|V(G)|, 2\delta\}$.
\end{thm}

Before discussing a 3-connected graph, we can show why the conjecture is not true for a 2-connected graph.  We construct a graph as follows by partitioning $V(G)$ into two parts $X$ and $Y$, as illustrated in Figure 7, where $G[X]$ is a cycle, $X = \{x_0, x_1, x_2\}$, and $Y$ is a path of length 3.  

As illustrated, the bond size of $[X,Y]$ is 10, which is a largest bond. Put 100 vertices in segments $[x_0, x_1]$,  $[x_1, x_2]$, and $[x_0, x_2]$,  to make the  longest cycle of the new graph. The largest bond size is still 10. Then, in this graph, the longest cycle does not meet the largest bond, so the conjecture does not hold for 2-connected graphs.

\vspace{-2mm}
\begin{center} 
\includegraphics[width=70mm, height=45mm]{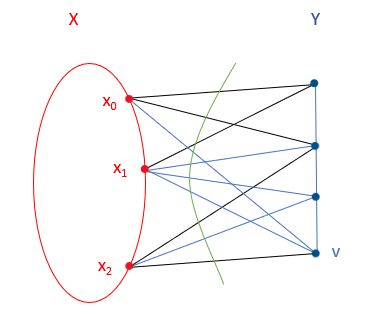}\\
Figure 7: A 2-connected graph
\end{center}

\noindent \textbf{\bf Proof of Theorem \ref{main}:} \\
\indent Suppose the theorem is $not$ true. Assume $G$ is a 3-connected graph with $n$ vertices and $m$ edges with a largest bond $B =[X, Y]$, and $X$ contains a longest cycle $C$, whose size is $c(G)$. Then $c(G)\le n$, $c^*(G)\le m - n + 2$, and $E(B)\cap E(C) = \emptyset$.

The proof will be divided into 8 cases: $c(G) = n, n - 1, n - 2, n - 3$, and $c^*(G) = m - n + 2, m - n + 1, m - n, m - n - 1$. 
\vspace{0.2in}

$\textbf{Case 1.1: $c(G) = n$}$ 

Since $c(G) = n$, the graph $G$ is a Hamiltonian graph with $B =[X, Y]$. In the Hamiltonian graph, the cycle $C$ must contain all the vertices. Therefore, $V(C) = X\cup Y$, thus $E(B)\cap E(C)\neq \emptyset$, a contradiction. 
\vspace{0.2in}

$\textbf{Case 1.2: $c(G) = n -1$}$ 

When $c(G) = n - 1$, we have the following proposition:

\textbf{Proposition}: Let $G$ be a 3-connected graph, and $B$ be a largest bond of $G$ with partition $(X, Y)$ of $V(G)$, then $|X|\ge 2$ and $|Y|\ge 2$. 

In Figure 8, suppose $|Y| = 1$, and $Y = \{u\}$. Let $uv\in E(G)$, and $Y_1$ = $Y\cup \{v\}$, and $X_1$ = $X - \{v\}$. Then $|[X_1, Y_1]|\ge |B| - 1$ + 2 = $|B|$ + 1. Clearly, since $G$ is 3-connected, $G[X_1]$ is connected, and $G[Y_1] $ is also connected. The new bond is a larger bond. This is a contradiction to the assumption that $B$ is a largest bond. Therefore, $|Y|\ge 2$. Similarly, $|X|\ge 2$.

\begin{center} 
\includegraphics[width=70mm, height=45mm]{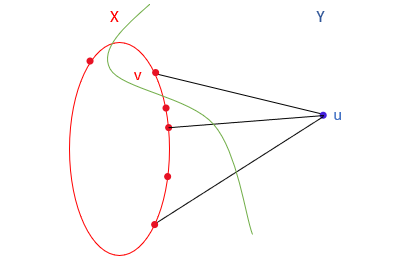}\\
Figure 8: A graph with an $(n - 1)$ cycle
\end{center}

$\textbf{Case 1.3: $c(G) = n -2$}$

Suppose $C$ is a longest cycle with length $n - 2$, and $B$ is a largest bond, such that $E(B)\cap E(C) = \emptyset$. Moreover, there is a partition $(X, Y)$ of $V(G)$, such that $V(C)\subseteq X$, and $B = [X, Y]$. According to the proposition above, $|X|, |Y|\ge 2$, thus when $c(G) = n - 2$, $|Y| = 2$, and $c(G) = |X|$.

In Figure 9, $Y = \{u, v\}$. As $G$ is 3-connected, each of $u, v$ has at least two neighbors in $X$. Choose two distinct vertices $a$,$d\in V(C)$, and they are as close as possible in $C$, such that $ua$, $vd\in E(G)$. We show next $ud$, $av\in E(G)$ too.

\textbf{Lemma 1}: $av\in E(G)$. By symmetry, $ud\in E(G)$.

Let $X_1 = X - \{a\}$, $Y_1 = Y\cup \{a\}$, if $va$ is not an edge, $|[X_1, Y_1]| $ will lose $ua$ but gain at least two more edges (note that $[X_1, Y_1]$ is a bond). This is a contradiction. Thus, the vertex $v$ must connect vertex $a$. Similarly, $ud\in E(G)$.

\begin{center} 
\includegraphics[width=70mm, height=45mm]{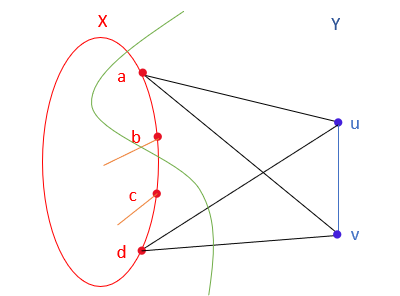}\\
Figure 9: A graph with an $(n - 2)$ cycle
\end{center}

Note that $ad$ is not an edge of $C$. Otherwise the longest cycle would include $u$ and $v$, to form a Hamiltonian cycle. Thus, there must be some vertices between $a$ and $d$. In fact, there must be at least two vertices $b, c$ between $a$ and $d$. Indeed, if there is only one vertex $b$ between $a$ and $d$, then replace the path $abd$ by the path $auvd$, and we obtain a longer cycle than $C$; a contradiction.

If $bu$ is an edge, we obtain a cycle longer than $C$ by replacing $bu$ with $aub$; a contradiction. Similarly, $bv$ cannot be an edge. Besides $a$ and $c$, vertex $b$ also connects to another vertex in $C$, as $d(b)\ge 3$. 

Now, let $Y_2 = \{u, v, a, b\}$ and $X_2 = X - \{a, b\}$, then clearly, both $G[X_2]$ and $G[Y_2]$ are connected, thus $[X_2, Y_2]$ is a bond. Moreover, $|[X_2, Y_2]|\ge  |B| - 2$ + 3 = $|B|$ + 1. This is a contradiction. 
\vspace{0.2in}

$\textbf{Case 1.4: $c(G) = n -3$}$ 

Suppose $C$ is a longest cycle with length $n - 3$, and $B$ is a largest bond, such that $E(B)\cap E(C) = \emptyset$. There is a partition $(X, Y)$ of $V(G)$, such that $V(C)\subseteq X$, and $B = [X, Y]$.  Then there are two possibilities for $Y$: $|Y|$ = 3 or $|Y|$ = 2.

\textbf{Case 1.4.1}: As in Figure 10, $|Y| = 3$, $Y = \{u, v, w\}$. As $G[Y]$ is connected, we may assume $uv, vw\in E(G)$, but $uw$ may be an edge or not. 

As $G$ is 3-connected, there are three disjoint paths between $Y$ and $V(C)$. Assume $uu_1$, $vv_1$, $ww_1\in E(G)$, where $u_1, v_1, w_1\in V(C)$ and $u_1, v_1, w_1$ are all distinct vertices. Consider the $[u_1, v_1]$-segment of $C$, $C[u_1, v_1]$, which does not contain $w_1$.

\begin{center} 
\includegraphics[width=65mm, height=45mm]{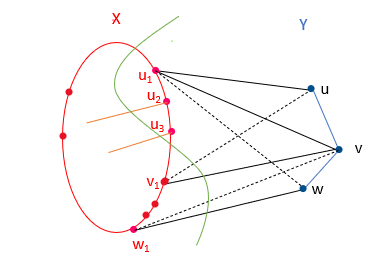}\\
Figure 10: A graph with an $(n - 3)$ cycle
\end{center}
\vspace{-2mm}

In $C[u_1, v_1]$, there must be at least two vertices $u_2$ and $u_3$ between $u_1$ and $v_1$. Otherwise we can find a longer cycle than $C$. Assume $u_1u_2, u_2u_3\in E(C)$. Note that $u_2$ cannot be adjacent to any of the vertices in $\{u, v, w\}$, and $u_3$ cannot be adjacent to $v$ or $w$. If $u_3 u\in E(G)$, then let $X_1 = X - \{u_3\}$, $Y_1 = Y \cup \{u_3\}$. Then $[X_1, Y_1]$ is a bond and $|[X_1, Y_1]|\ge  |B| - 1$ + 2 = $|B|$ + 1; a contradiction. Therefore, $u_3u\not\in E(G)$ as well. Note that $d(u_2), d(u_3)\ge 3$. If $u_1u_3\in E(G)$, let  $X_1 = X - \{u_1, u_2\}$, and $Y_1 = Y\cup \{u_1, u_2\}$. If $u_1u_3\not\in E(G)$, let $Y_1 = \{u, v, w, u_1, u_2, u_3\}$ and $X_1 = X - \{u_1, u_2, u_3\}$, then in either case, we get a new bond $[X_1, Y_1]$, so that  $|[X_1, Y_1]|\ge |B| - 3$ + 4 = $|B|$ + 1. This is a contradiction.

\textbf{Case 1.4.2}: $|Y| = 2$, $|X| = n - 2$. Suppose $X = V(C)\cup \{w\}$, $Y = \{u, v\}$, $B = [X, Y]$. Clearly, $uv\in E(G)$ as $V(Y)$ is connected. As $d(u), d(v)\ge 3$, each of $u$ and $v$ has at least two neighbors in $X$.

As shown in Figure 11, suppose vertex $w$ has exactly one neighbor in $Y$, say $wv\in E(G)$, and thus $w$ has at least two neighbors in $C$, thus $wa, wb\in E(G)$ for some $a, b\in V(C)$ . Then the bond $[X - \{w\}$, $Y\cup \{w\}]$ has at least one more edge than $[X, Y]$. This is a contradiction. Therefore, $wu, wv\in E(G)$, or $w$ is not adjacent to any vertex in $Y$ at all.

\begin{center} 
\includegraphics[width=65mm]{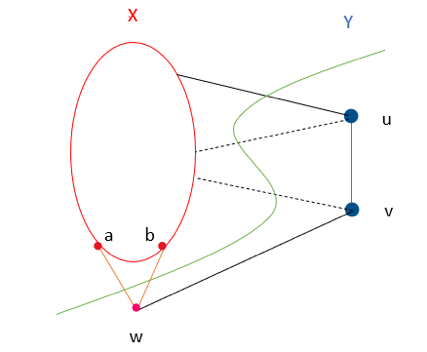}\\
Figure 11: A graph with an $(n - 3)$ cycle
\end{center}

First, assume $w$ has no neighbors in $Y$, as shown in Figure 12. Thus, all neighbors of $w$ are in $X$. As $G$ is 3-connected, there are two distinct vertices $u_1, v_1\in V(C)$, such that $uu_1, vv_1\in E(G)$. By Lemma 1, $u_1v, v_1u\in E(G)$. Moreover, at least two vertices $u_3$ and $u_4$ exist in between $u_1$ and $v_1$ in $C$. Here $u_1u_3\in E(G)$ and $u_3$ is not adjacent to any of $u$ and $v$; otherwise we get a cycle longer than $C$. As $d(u_3)\ge 3$, $u_3$ is adjacent to some vertices in $V(G) - \{u, v, u_1, u_4\}$.

Let $X_1 = X - \{u_1, u_3\}$, $Y_1 = Y\cup \{u_1, u_3\}$. Then $[X_1, Y_1]$ is a bond with size at least $|B| - 2 + 3 = |B|$ + 1. This is a contradiction again.

\begin{center} 
\includegraphics[width=65mm, height=42mm]{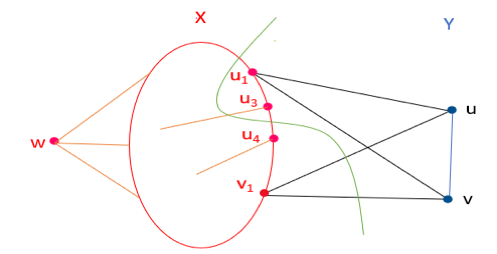}\\
Figure 12: A graph with an $(n - 3)$ cycle
\end{center}

Finally, as shown in Figure 13, $w$ has at least one neighbor in $C$ and two neighbors in $Y$. As $G$ is 3-connected, $\{u, v\}$ has at least two neighbors $u_1, v_1$ in $C$. Moreover, if $w$ has a unique neighbor $a$ in $C$, then we can choose $u_1, v_1$ different from $a$. Otherwise $G$ would have a vertex-cut of size 2; a contradiction as $G$ is 3-connected. Now similar to Lemma 1, we can show $uu_1, uv_1, vu_1, vv_1$ are all edges of $G$. There are at least two vertices $u_2$, $u_3$ between $u_1$ and $v_1$ in $C$, where $u_1u_2\in E(G)$, otherwise $G$ has a cycle longer than $C$. Moreover, $u_2u, u_2v\not\in E(G)$. Otherwise $G$ has a cycle longer than $C$ again. 

Let $Y_1 = \{u, v, u_1, u_2\}$ and $X_1 = X - \{u_1, u_2\}$, then $|[X_1, Y_1]|\ge |B| - 2 + 3 = |B|$ + 1. This is a contradiction.

\begin{center} 
\includegraphics[width=75mm, height=45mm]{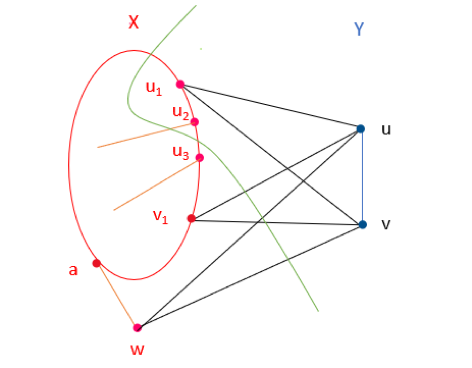}\\
Figure 13: A graph with an $(n - 3)$ cycle
\end{center}

Therefore, if $c(G)\ge n - 3$, the longest cycle $C$ and largest bond $B$ must meet. 
\vspace{0.2in}

$\textbf{Case 2.1: $c^*(G) = m - n + 2$}$ 

Note that $G$ is dual Hamiltonian. $B$ has size $|E(G)| - |V(G)|$ + 2. Since it is a dual Hamiltonian graph, $G[X]$ and $G[Y]$ are both trees. Thus, any cycle must meet $B$; a contradiction.

$\textbf{Case 2.2: $c^*(G) = m - n + 1$}$ 

Assume $E(B)\cap E(C) = \emptyset$ and without loss of generality, assume $V(C)\subseteq X$, then $G[X]\cup G[Y]$ has $m - (m - n + 1) = n - 1 = |X| + |Y| - 1 = |X| + (|Y| - 1)$ edges.

$G[Y]$ has at least $|Y| - 1$ edges, but $G[X]$ has at least $|X|$ edges since $X$ contains the cycle $C$. Thus, $G[X]$ has exactly $|X|$ edges, and $G[X]$ has exactly one cycle, $C$, and $G[Y]$ is a tree.

As shown in Figure 14, if we delete an edge $e$ of cycle $C$ from $G[X], G[X] - e$ will be a tree. By Theory $\ref{dirc}$, then $|C|\ge 2\delta\ge 6$. We take $u$ and $v$, two adjacent vertices in $C - e$. Since $G$ is a 3-connected graph in $G - e$, $d(u), d(v)\ge 3$. Thus, there are $a, b\in V(G)\backslash V(C)$, neither of which is incident to $e$, such that $ua, vb\in E(G)$. Note $a\not= b$ as $C$ is the largest cycle of $G$. If $a\in X$, traverse $u$ in the tree $G[X] - e$ along $a$ to reach a pendant $a_1$. Similarly, if $b\in X$, traverse $v$ in the tree $G[X] - e$ along $b$ to reach a pendant $b_1$. As $d(a_1), d(b_1)\ge 3$,  $a_1$ and $b_1$ are both connected to some vertices of $Y$. As $G[Y] $ is connected, a new longer cycle forms by replacing $uv$ with a path from $u$ to $a_1$, to $Y$, back to $b_1$, then to $v$ as highlighted in blue. This is a contradiction again. If $a\in Y$ or $b\in Y$, one gets a contradiction similarly. 

\vspace{-2mm}

\begin{center} 
\includegraphics[width=75mm, height=55mm]{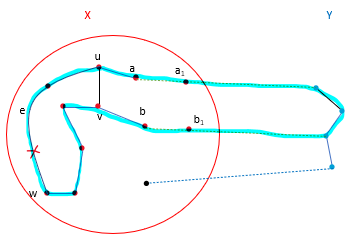}\\
Figure 14: A graph with an $(m - n + 1)$-size bond
\end{center}

$\textbf{Case 2.3: $c^*(G) = m - n$}$ 

As shown in Figure 15, $G[X]\cup G[Y]$ has exactly $|X| + |Y|$ edges. If $G[Y]$ has a cycle, then $G[X]$ has exactly $|X|$ edges, and the proof is similar to that of Case 2.2. Thus, we assume that $G[Y]$ is a tree with $|Y| - 1$ edges, then $G[X]$ has $|X| + 1$ edges, which means $G[X]$ is a tree plus two edges. Remove an edge f from $C$, and another edge $e$ not incident to $f$, such that $G[X] - e - f$ is a tree.

\begin{center} 
\includegraphics[width=85mm]{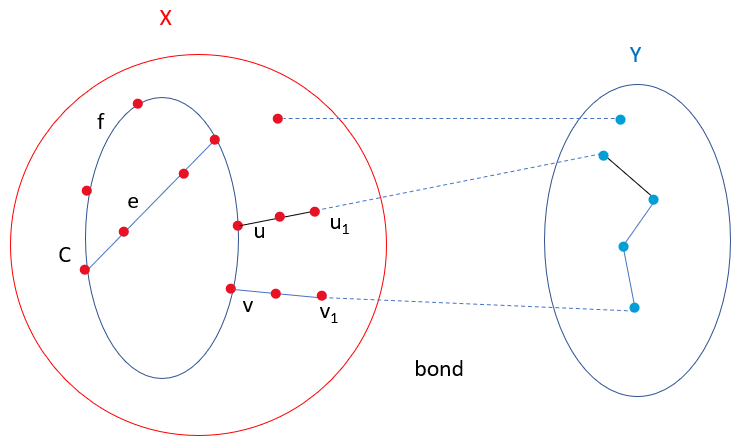}\\
Figure 15: A graph with an $(m - n)$-size bond
\end{center}

Now $G - e - f$ has minimum degree at least 2 as $G$ is 3-connected. Take two adjacent vertices $u, v$ in the path $C - e$, such that the degrees of both vertices are at least 3 in $G - e - f$. Note that such vertices must exist as $|C|\ge 6$. Consider the tree $G[X] - e - f$, similar to the last case, traverse $u$ and $v$ outside of $C$ to two pendants $u_1$ and $v_1$ in the tree if $u$ is not adjacent to a vertex in $Y$. As the degree of $u_1$ and $v_1$ are at least 2 in $G - e - f$, $u_1$ has a neighbor in $Y$, so does $v_1$. However, since $G[Y]$ is connected, there is a path connecting these two vertices. Now we easily get a longer cycle than $C$; a contradiction. If $u$ or $v$ is adjacent to a vertex in $Y$, one gets a longer cycle than $C$ again using a similar argument; a contradiction.

$\textbf{Case 2.4: $c^*(G) = m - n - 1$}$ 

In this case, $G[X]\cup G[Y]$ contains $|X| + |Y| + 1$ edges. If $G[X]$ has at most $|X| + 1$ edges, then we get a contradiction as in Case 2.2 and Case 2.3. Thus, without loss of generality, assume $G[Y]$ is a tree with $|Y| - 1$ edges, then $G[X]$ has $|X| + 2$ edges, which means $G[X]$ is a tree plus 3 edges. 

A chordal path of $C$ is a path starting and ending with vertices of $C$ and internally disjoint from $C$. There are two cases to be discussed. 

\textbf{Case 2.4.1}: $C$ has two chordal paths, as illustrated in Figure 16. Note that again, by Dirac’s Theorem \ref{dirc}, $c(G)\ge 6$ as $G$ is 3-connected.

\begin{center} 
\includegraphics[width=78mm, height=43mm]{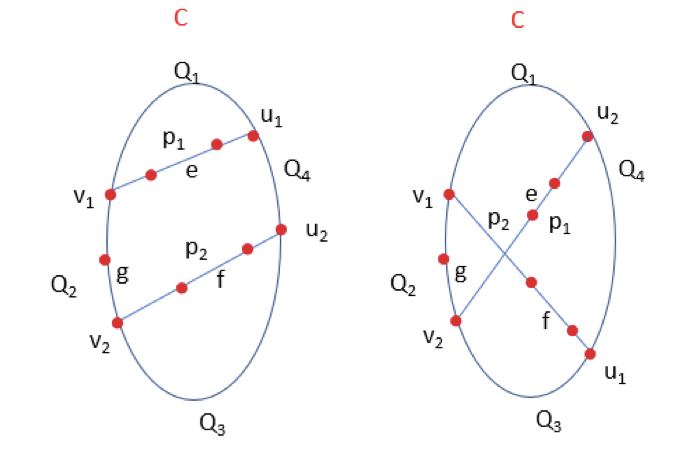}\\
Figure 16: A cycle with two chordal paths 
\end{center}

These two chordal paths divide cycle $C$ into four paths $Q_1, Q_2, Q_3$ and $Q_4$ (some of the paths can be empty). We assume $Q_1$ is the shortest path. Remove edges $e$ and $f$ from the two chordal paths, respectively, and an edge $g$ from the shortest path with at least one edge. Then $G - e - f - g$ is a tree. 

\begin{center} 
\includegraphics[width=70mm]{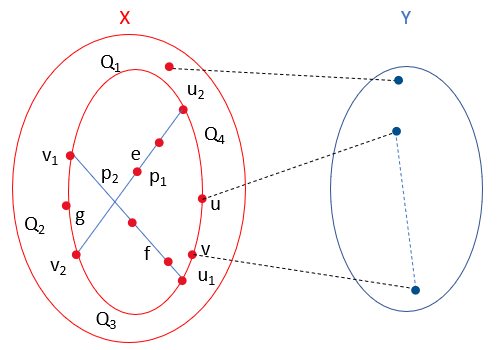}\\
Figure 17: A graph with two chordal paths
\end{center}
\vspace{-2mm}

If $C - g$ has two adjacent vertices $u$ and $v$ with degree at least 3 in $G - e - f - g$, as shown in Figure 17, we proceed to traverse from $u$ and $v$ to two pendants; a contradiction as in the previous cases. Thus, we may assume we cannot find two such vertices. Therefore, each of the $Q_i (1\le  i\le 4)$  has at most one internal vertex, and each chordal path is an edge with both end vertices having degree 3 (when $e$ and $f$ are not incident). Moreover, if $e, f$ are incident to a vertex $t$ in $C$, then $d(t)$ = 4. As $G$ is 3-connected, at least three paths, say $Q_2, Q_3$ and $Q_4$ all have exactly one internal vertex. Thus, we can find two interval vertices $a$ and $b$ of $V(C)$ in $Q_4$ and $Q_3$, respectively, as shown in Figure 18. As $G - e - f - g$ is a tree, $a, b$ is not adjacent to any vertex in $V(C)$.

\begin{center} 
\includegraphics[width=82mm]{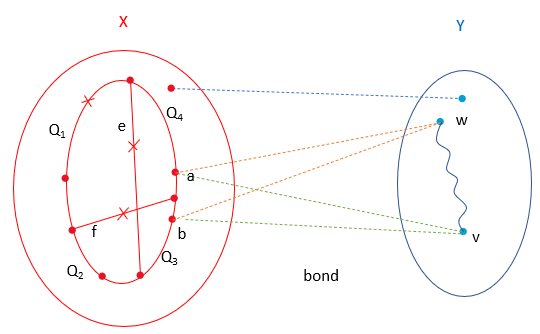}\\
Figure 18: A graph with an $(m - n - 1)$-size bond
\end{center}

Note that each $(a, b)$-path disjoint from $C$ must have length at most 2. As $d(a), d(b)\ge 3$, there are $a_1, b_1\in V(G)$ such that $aa_1, bb_1\in E(G)$. As in the previous case, we can find an $(a, b)$-path internally disjoint from $C$. Therefore, $a, b$ must be adjacent to a common vertex $v$ in $Y$; otherwise $G$ has a longer cycle than $C$.  If $d(a)$ = 3, then $[X - \{a\}$, $Y\cup \{a\}]$ is a larger bond than $[X, Y]$; a contradiction. Thus $d(a)\ge 4$. Similarly, $d(b)\ge 4$. Therefore, $a, b$ must be adjacent to another vertex $w$ in $Y$ by the above argument.

However, since $G[Y]$ is connected, there is a $(v, w)$-path in $G[Y]$. Now we can easily find an $(a, b)$-path internally disjoint from $C$ and longer than 2; a contradiction.

\textbf{Case 2.4.2}: $C$ has at most one chordal path as two scenarios, as shown in Figure 19. We can delete three edges $e, f, g$ from $G[X]$ to form a tree again, where $e$ is an edge of the possible chordal path, $g$ is an edge of $C$, and $f$ is an edge of another cycle in $G[X]$.

\indent In both scenarios, we can find two adjacent vertices $a, b$ of $V(C)$ with degree at least 3 in $G - e - f - g$. Proceed as before, we can find a longer cycle than $C$. This is a contradiction. 

This completes the proof of the final case, and thus the proof of Theorem \ref{main}.

\begin{center} 
\includegraphics[width=82mm]{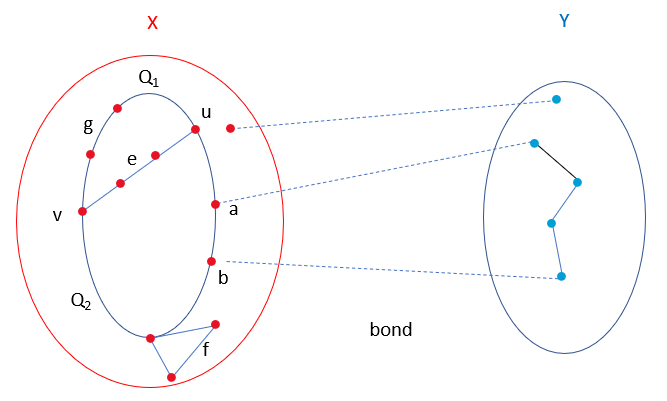}\\
\includegraphics[width=82mm]{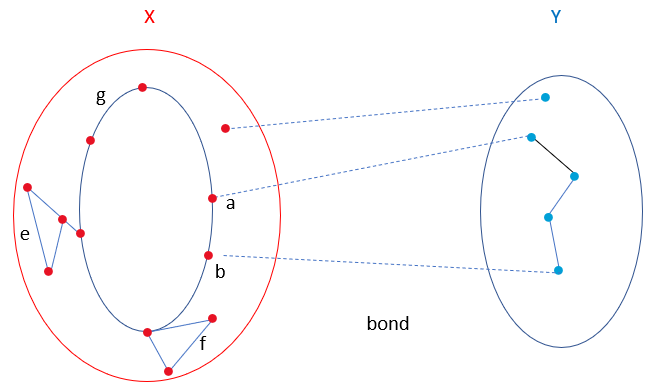}\\
Figure 19: A graph with an $(m - n - 1)$-size bond
\end{center}

\section{Conclusion and Further Study}

Research in cycles and related topics is a fundamental area in graph theory, and there is much interest in studying the longest cycles and largest bonds in graphs.  In this research, we have partially proven a conjecture involving the longest cycles and largest bonds in 3-connected graphs. We plan to continue to investigate the general case. We have also found the co-spectrum of the generalized Petersen graph $P(n, k)$. 

Since many real-world problems can be represented by graphs, graph theory has found extensive application in biochemistry, bioinformatics, computer science and artificial intelligence, social science, accounting, business, and many other diverse areas. 

The current research is one step forward in proving the conjecture that in a 3-connected graph, any longest cycle must meet any largest bond. A lot of further mathematical research can be conducted following the current research, so that we can finally prove the conjecture true in the general case.

\section*{Acknowledgements}

I would like to convey my deepest appreciation to Dr. Haidong Wu, professor of mathematics at the University of Mississippi, for his great support and guidance. He has been inspiring me to pursue my interest in Discrete Math and Graph Theory since May 2021. A big thank-you to all my teachers and schoolmates at Diamond Bar High School for their support. My heartfelt appreciation also goes to my parents for always being there and supporting me to take on new challenges. This research project would not be possible without you all. Thank you!

\end{document}